\documentclass[11pt]{article}
\usepackage[latin1]{inputenc}
\usepackage{amsmath,amssymb,indentfirst}
\usepackage{graphicx, epsfig}
\newcommand {\tr} {\mbox{tr}} 
\newcommand{\bfsym}{\boldsymbol}  
     
\newcommand {\J}{\widetilde{\!\!\widetilde{J}}}

\def \diag {\mbox{diag}}

\title{Adjusted likelihood inference in an elliptical multivariate errors-in-variables model}
\author
{
   Tatiane F. N. Melo\\
     {\it \footnotesize Instituto de Matemática e Estatística,
     Universidade Federal de Goiás, Brazil}
     \vspace{-0.2cm}\\
\\ \\
   Silvia L. P. Ferrari\footnote{Corresponding author. Email: {\tt silviaferrari.usp@gmail.com}}\\
     {\it \footnotesize Departamento de Estatística,
     Universidade de São Paulo, Brazil}
     \vspace{-0.2cm}\\
}
\date{}

\begin{document}
\maketitle

{\footnotesize 
\noindent{\bf Abstract:} 

\noindent
In this paper we obtain an adjusted version of the likelihood 
ratio test for errors-in-variables multivariate linear regression 
models. The error terms are allowed to follow a multivariate 
distribution in the class of the elliptical distributions, 
which has the multivariate normal distribution as a special case. 
We derive a modified likelihood ratio statistic that follows a 
chi-squared distribution with a high degree of accuracy. 
Our results generalize those in Melo \& Ferrari~
({\it Advances in Statistical Analysis\/}, 2010, {\bf 94}, 
75--87) by allowing the parameter of interest to be 
vector-valued in the multivariate errors-in-variables model. 
We report a simulation study which shows that the proposed 
test displays superior finite sample behavior relative to the 
standard likelihood ratio test. 

\vspace{0.1cm}
}
\section{Introduction}\label{sec1}

Statisticians are often faced with the problem of modeling data measured with error. 
As an example, we refer to Aoki et al.~(2001), who compared the effectiveness of two types of 
toothbrushes in removing dental plaque. One explanatory variable is the dental plaque index before 
toothbrushing and the response variable is the dental plaque index after brushing, the amount of 
plaque being imprecisely measured. The authors proposed a null intercept regression model that 
assumes that this explanatory variable is measured with an additive random error and the measurement 
error of the response variable is assumed to be absorbed by the error term of the model.

Errors-in-variables models are generalizations of classical regression models. The true (non-observable) explanatory variables are treated either as random variables, in which case the model is said to be 
structural, or as unknown parameters, leading to a functional model. Structural models are, in general, non-identifiable, while functional models induce unlimited likelihood functions. 
Such difficulties disappear if some variances are assumed to be known (e.g. Chan \& Mak~(1979) and 
Wong~(1989)) or that the intercept is null (e.g. Aoki et al.~(2001)).
For details on errors-in-variables models the reader is referred to, for instance, Fuller~(1987) and 
Buonaccorsi~(2010). 

The most popular errors-in-variables models for continuous outcomes are based on normality assumptions.
The family of the elliptical distributions provides a useful alternative to the normal distribution
when outlying observations are present in the data. It nests the normal distribution, 
heavy-tailed distributions, such as the exponential power and the Student-t distributions, and
light-tailed distributions. Further information on elliptical distributions can be found in Fang et al. (1990) 
and Fang and Anderson (1990).

As shown by Melo \& Ferrari~(2010), statistical inference in errors-in-variables models
based on first-order asymptotic approximations can be imprecise for small or 
moderate sized samples. In particular the type I error of the likelihood ratio test 
is often larger than the designed level of the test. Skovgaard~(2001) proposed a 
general strategy to adjust the likelihood ratio statistic
when interest lies in inference on a vector-valued parameter. The adjustment makes the 
resulting statistic to follow a chi-squared distribution with a high degree of accuracy. The adjustment is  
broadly general, but requires either some unusual likelihood quantities or the identification of a
suitable ancillary statistic such that, when coupled with the maximum likelihood estimator, constitutes 
a sufficient statistic for the model. In the present paper, we obtain an appropriate ancilary statistic 
and derive Skovgaard's adjustment for a structural elliptical multivariate errors-in-variables model.  

The paper unfolds as follows. Section \ref{sec2} introduces the model. Section \ref{sec3} 
contains our main results, namely the ancillary statistic and an explicit formula for the modified 
likelihood ratio test. The finite 
sample behavior of the likelihood ratio test and its adjusted version is evaluated and discussed
in  Section \ref{sec4}. Our simulation results clearly show that the likelihood ratio test 
tends to be oversized and its modified version is much less size-distorted. 
Finally, Section \ref{sec6} closes the paper with our conclusions. 
Technical details are left for three appendices. 

\section{The model}
\label{sec2}

The $(l+1) \times 1$ random vector ${\bfsym Z}$ is said to have
a $(l+1)$-variate elliptical distribution with location vector
$\bfsym{\mu}$ $((l+1) \times 1)$, dispersion matrix $\Sigma$ $((l+1) \times (l+1))$
and density generating function $p_0$, and we write
${\bfsym Z} \sim El_{(l+1)}(\bfsym{\mu}, \Sigma; p_0)$, if
\begin{eqnarray*}\label{E.3.2.4}
{\bfsym Z} \stackrel{d}{=} \bfsym{\mu} + A {\bfsym Z}^*,
\end{eqnarray*}
where $A$ is  $(l+1) \times k$ matrix with ${\rm rank}(A) = k$, $A A^\top = \Sigma$ 
and ${\bfsym Z}^*$ is a $(l+1) \times 1$ random vector with density
function $p_0({\bfsym z}^\top {\bfsym z})$, for ${\bfsym z} \in {\Re}^{(l+1)}$.
The notation ${\bfsym X} \stackrel{d}{=} {\bfsym Y}$ indicates that ${\bfsym X}$ 
and ${\bfsym Y}$ have the same distribution.
It is assumed that
$\int_{0}^{\infty} y^{(l+1)/2-1} p_0(y) dy < \infty.$
The density function of ${\bfsym Z}$ is
\begin{eqnarray}\label{E.3.2.7}
p({\bfsym z},\mu,\Sigma) = |\Sigma|^{-1/2} p_0\left(({\bfsym z} 
- \bfsym{\mu})^\top \Sigma^{-1} ({\bfsym z} - \bfsym{\mu})\right).
\end{eqnarray}
Some special cases of (\ref{E.3.2.7}) are the following multivariate
distributions: normal, exponential power, Pearson II, Pearson VII,
Student-$t$, generalized Student-$t$, logistic I, logistic II and Cauchy.
The elliptical distributions share many properties with the
multivariate normal distribution. In particular,
marginal distributions are elliptical.
For a full account of the properties of the elliptical distributions, see
Fang et al.~(1990, Sect. 2.5).

We consider the following model, which consists of $p$ independent errors-in-variables structural models:
\begin{equation}\label{E.5.1}
\begin{split}
{\bfsym Y}_{jk} &= {\bfsym \alpha}_k + \bfsym{\beta}_k x_{jk} + {\bfsym e}_{jk}, 
\\
X_{jk} &= x_{jk} + u_{jk},
\end{split}
\end{equation}
for $j = 1, 2, \ldots, n_k$ and $k = 1, 2, \ldots, p$, where ${\bfsym Y}_{jk} = \left(Y_{1jk}, Y_{2jk}, 
\ldots, Y_{ljk}\right)$, ${\bfsym \alpha}_k = \left(\alpha_{1k}, \alpha_{2k}, \ldots, \alpha_{lk}\right)$, 
$\bfsym{\beta}_k = \left(\beta_{1k}, \beta_{2k}, \ldots, \beta_{lk}\right)$ and ${\bfsym e}_{jk} = 
\left(e_{1jk}, e_{2jk}, \ldots, e_{ljk}\right)^\top$. Here,  $x_{jk}$ is not observed directly.
Instead, we observe $X_{jk}$, which is viewed as $x_{jk}$ plus a measurement error, $u_{jk}$. 
We assume that ${\rm E}({\bfsym e}_{jk}) = {\bfsym 0}$, ${\rm Var}
({\bfsym e}_{jk}) = \Sigma_{e_k} = \diag\left\{\sigma^2_{e_{1k}}, \sigma^2_{e_{2k}}, \ldots, 
\sigma^2_{e_{lk}}\right\}$ and ${\rm Cov}({\bfsym e}_{j'k'}, {\bfsym e}_{jk}) = 0$; ${\rm E}(u_{jk}) = 0$, 
${\rm Var}(u_{jk}) = \sigma^2_{u_k}$ and ${\rm Cov}(u_{j'k'}, u_{jk}) = 0$; ${\rm E}(x_{jk}) = \mu_{x_k}$, 
${\rm Var}(x_{jk}) = \sigma^2_{x_k}$ and ${\rm Cov}(x_{j'k'}, x_{jk}) = 0$; ${\rm Cov}(x_{j'k'}, u_{jk}) = 0$; 
${\rm Cov}({\bfsym e}_{j'k'}, u_{jk}) = {\rm Cov}({\bfsym e}_{j'k'}, x_{jk}) = {\bfsym 0}$, with $(j', k') 
\ne (j, k)$.

Model (\ref{E.5.1}) can be written as
\begin{equation}\label{E.5.3}
{\bfsym Z}_{jk} = {\bfsym\delta}_k + \Delta_k {\bfsym b}_{jk},
\end{equation}
for $j = 1, 2, \ldots, n_k$ and $k = 1, 2, \ldots, p$, where
\begin{equation*}\label{E.5.4}
{\bfsym Z}_{jk} = \left(\begin{array}{c} {\bfsym Y}_{jk} \\ X_{jk} \end{array}\right), \:\:\: 
{\bfsym\delta}_k = \left(\begin{array}{c} {\bfsym \alpha}_k \\ 0 \end{array}\right), \:\:\: 
\Delta_k = \left(\begin{array}{ccc} \bfsym{\beta}_k&I_l&{\bfsym 0} \\ 1&{\bfsym 0}&1 
\end{array}\right) \:\:\: {\rm and}\:\:\: {\bfsym b}_{jk} = \left(\begin{array}{c} x_{jk} \\ 
{\bfsym e}_{jk} \\ u_{jk} \end{array}\right),
\end{equation*}
where $I_l$ is the identity matrix of dimension $l$. We assume that, for each $k = 1,2, \ldots, p$, 
the errors ${\bfsym b}_{1k}, {\bfsym b}_{2k}, \ldots, {\bfsym b}_{{n_k}k}$ are independent and
${\bfsym b}_{jk} \sim El_{(l+2)}(\bfsym{\eta}_k, \Omega_k; p_0)$, with
\begin{eqnarray*}\label{E.5.5}
{\bfsym \eta}_k = \left(\begin{array}{c}
\mu_{x_k} \\ {\bfsym 0} \\ 0
\end{array}\right)
\:\:\:{\rm and}\:\:\:
\Omega_k = \left(\begin{array}{ccc} 
\sigma^2_{x_k} & {\bfsym 0} & 0 \\ {\bfsym 0} & \Sigma_{e_k} & {\bfsym 0} \\ 0 & {\bfsym 0} & \sigma^2_{u_k}
\end{array}\right).
\end{eqnarray*}
Therefore, for each $k = 1,2, \ldots, p$, the random vectors ${\bfsym Z}_{1k}, {\bfsym Z}_{2k}, 
\ldots, {\bfsym Z}_{{n_k}k}$ are independent and ${\bfsym Z}_{jk} \sim El_{l+1}({\bfsym \mu}_k, 
\Sigma_k; p_0)$, with ${\bfsym \mu}_k = {\bfsym \delta}_k + \Delta_k {\bfsym \eta}_k$ and 
$\Sigma_k = \Delta_k \Omega_k \Delta_k^\top$ (Fang et al. 1990, Sect. 2.5). We can write 
${\bfsym \mu}_k$ and $\Sigma_k$ as
\begin{eqnarray*}\label{E.5.6}
{{\bfsym \mu}_k} = \left(\begin{array}{c}
{\bfsym \alpha}_k + {\bfsym \beta}_k \mu_{x_k} \\ \mu_{x_k} 
\end{array}\right)
\:\:\:{\rm and}\:\:\:
\Sigma_k = \left(\begin{array}{cc} 
{\bfsym \beta}_k \sigma^2_{x_k} {\bfsym \beta}_k^\top + \Sigma_{e_k} & {\bfsym \beta}_k \sigma^2_{x_k} \\ 
\sigma^2_{x_k} {\bfsym \beta}_k^\top & \sigma^2_{x_k} + \sigma^2_{u_k}
\end{array}\right).
\end{eqnarray*}

Regression model (\ref{E.5.1}) generalizes the normal structural models proposed by Cox~(1976)
($l=1$) and Russo et al.~(2009) ($l=2$), and the elliptical structural model considered in Melo \& Ferrari~(2010) ($l = p = 1$); a closely related model is presented by Garcia-Alfaro \& Bolfarine~(2001). 
As expected, model (\ref{E.5.1}) is not identifiable because the relation between 
the parameters of the distribution of ${\bfsym Z}_{jk}$ and ${\bfsym \theta}_{(k)} = 
\big({\bfsym \beta}_k^\top, \alpha_k, \mu_{x_k}, \sigma^2_{x_k}, \sigma^2_{u_k}, 
{\bfsym \sigma}^{2^\top}_{e_k} \big)^\top$ is not unique. Assumptions on 
$\sigma^2_{x_k}$ and ${\bfsym \sigma}^2_{e_k}$ are usually imposed to overcome 
identifiability problems. It is common to assume that the $\lambda_{x_k} = 
\sigma^2_{x_k}/\sigma^2_{u_k}$ or ${\lambda}_{e_{ik}} = {\sigma}^2_{e_{ik}}/\sigma^2_{u_k}$, 
for $i = 1, 2, \ldots, l$, is known. An alternative assumption is that 
the intercept $\alpha_k$ is known (see Aoki et al., 2001). 
Under each of these identifiability assumptions we have:
\begin{description}
	\item[(i)]
   {if $\lambda_{x_k}$ is known,
    $$
    {\bfsym \theta}_{(k)} = \left({\bfsym \beta}_k^\top, \alpha_k, \mu_{x_k}, \sigma^2_{u_k}, 
    {\bfsym \sigma}^{2^\top}_{e_k}\right)^\top
    \ \ {\mbox{and}}\ \ 
    \Sigma_k = 
    \left(
    \begin{array}{cc} 
    {\bfsym \beta}_k {\bfsym \beta}_k^\top \lambda_{x_k}\sigma^2_{u_k} + \Sigma_{e_k} & 
    {\bfsym \beta}_k \lambda_{x_k}\sigma^2_{u_k} \\ 
    \lambda_{x_k}\sigma^2_{u_k} {\bfsym \beta}_k^\top & (\lambda_{x_k} + 1)\sigma^2_{u_k}
    \end{array}
    \right);
    $$
    }
    \item[(ii)]
    {if $\lambda_{e_k}$ is known,
    $$
    {\bfsym \theta}_{(k)} = \left({\bfsym \beta}_k^\top, \alpha_k, \mu_{x_k}, \sigma^2_{x_k}, 
    \sigma^2_{u_k}\right)^\top
    \ \ {\mbox{and}}\ \ 
    \Sigma_k = 
    \left(
    \begin{array}{cc} 
    {\bfsym \beta}_k \sigma^2_{x_k} {\bfsym \beta}_k^\top + {\lambda}_{e_k}\sigma^2_{u_k} & 
    {\bfsym \beta}_k \sigma^2_{x_k} \\ 
         \sigma^2_{x_k} {\bfsym \beta}_k^\top & \sigma^2_{x_k} + \sigma^2_{u_k}
    \end{array}
    \right),
    $$
    with ${\lambda}_{e_k} = \diag\left\{\lambda_{e_{1k}}, \lambda_{e_{2k}}, \ldots, 
    \lambda_{e_{lk}}\right\}$;
   }
   \item[(iii)]
    {if $\alpha$ is known,
    $$
    {\bfsym \theta}_{(k)} = \left({\bfsym \beta}_k^\top, \mu_{x_k}, \sigma^2_{x_k}, \sigma^2_{u_k}, 
    {\bfsym \sigma}^{2^\top}_{e_k}\right)^\top
    \ \ {\mbox{and}}\ \ 
    \Sigma_k = 
    \left(
    \begin{array}{cc} 
    {\bfsym \beta}_k \sigma^2_{x_k} {\bfsym \beta}_k^\top + \Sigma_{e_k} & {\bfsym \beta}_k 
    \sigma^2_{x_k} \\ 
    \sigma^2_{x_k} {\bfsym \beta}_k^\top & \sigma^2_{x_k} + \sigma^2_{u_k}
    \end{array}
    \right).
    $$
   }
\end{description}

The independent structural elliptical model can be defined in terms of the density function of 
${\bfsym Z} = \left({\bfsym Z}_{(1)}^\top, {\bfsym Z}_{(2)}^\top, \ldots, {\bfsym Z}_{(p)}^\top
\right)^\top$, with ${\bfsym Z}_{(k)} = \left({\bfsym Z}_{1k}^\top, {\bfsym Z}_{2k}^\top, \ldots, 
{\bfsym Z}_{{n_k}k}^\top\right)^\top$, for $k = 1, 2, \ldots, p$, which is given by
\begin{equation*}\label{E.5.9}
p_{\bfsym Z}({\bfsym z}, {\bfsym \theta}) = \prod_{k=1}^{p} \prod_{j=1}^{n_k} |\Sigma_k|^{-1/2} 
p_0({\bfsym d}_{jk}^\top \Sigma_k^{-1} {\bfsym d}_{jk}),
\end{equation*}
where ${\bfsym d}_{jk} = {\bfsym d}_{jk}({\bfsym \theta}_{(k)}) = {\bfsym z}_{jk} - {\bfsym \mu}_k$, 
for $j = 1, 2, \ldots, n_k$, $k = 1, 2, \ldots, p$, and ${\bfsym \theta} = \left({\bfsym 
\theta}_{(1)}^\top, {\bfsym \theta}_{(2)}^\top, \ldots, {\bfsym \theta}_{(p)}^\top\right)^\top$. 

The log-likelihood function for the $k$-th group, $k = 1, 2, \ldots, p$, is given by
\begin{eqnarray*}\label{E.5.10}
\ell_k({\bfsym \theta}, {\bfsym z}) = -\frac{n_k}{2}\log |\Sigma_k| + \sum_{j=1}^{n_k} \log p_0
({\bfsym d}_{jk}^\top \Sigma_k^{-1} {\bfsym d}_{jk}).
\end{eqnarray*}
For a sample of size  $n = \sum_{k=1}^{p} n_k$ and $p$ populations, the log-likelihood function is
\begin{eqnarray}\label{E.5.11}
\ell({\bfsym \theta}, {\bfsym z}) = \sum_{k=1}^{p} \ell_k({\bfsym \theta}, {\bfsym z}).
\end{eqnarray}
Maximum likelihood estimation of the parameters  can be carried out by numerically
maximizing the log-likelihood function (\ref{E.5.11}) through an iterative algorithm such as 
the Newton--Raphson, the Fisher scoring, EM or BFGS. Our numerical results were obtained using 
the library function {\tt MaxBFGS} in the {\tt Ox} matrix programming language (Doornik 2006). 

\section{Ancillary statistic and modified likelihood ratio test}
\label{sec3}

The parameter vector ${\bfsym \theta}$ is partitioned as ${\bfsym \theta} = ({\bfsym \psi}^\top, 
{\bfsym \omega}^\top)^\top$, with ${\bfsym \psi}$ representing the parameter of interest and ${\bfsym 
\omega}$ the nuisance parameter. Our interest lies in testing ${\cal H}_{0}: 
{{\bfsym \psi}} = {{\bfsym \psi}}^{(0)} \:\:\: {\rm versus} \:\:\: {\cal H}_{1}: {{\bfsym \psi}} 
\neq {{\bfsym \psi}}^{(0)}$, where ${{\bfsym \psi}}^{(0)}$ is a $q$-dimensional parameter of known
constants. The maximum likelihood estimator of ${\bfsym \theta}$ is denoted by
$\bfsym{\widehat{\theta}} = (\bfsym{\widehat{\psi}}^\top, \bfsym{\widehat{\omega}}^\top)^\top$ 
and the corresponding estimator obtained under the null hypothesis is  $\bfsym{\widetilde\theta} = 
(\bfsym{\widetilde\psi}^\top, \bfsym{\widetilde\omega} ^\top)^\top$, where $\bfsym{\widetilde\psi} = 
{\bfsym{\psi}}^{(0)}$. We use hat and tilde to indicate evaluation at $\bfsym{\widehat\theta}$ 
and $\bfsym{\widetilde\theta}$, respectively.

The likelihood ratio statistic for testing ${\cal H}_0$ is given by
\begin{eqnarray*}\label{E.4.11}
LR = 2\:\left\{\ell(\bfsym{\widehat{\psi}}) - \ell(\bfsym{\widetilde{\psi}})\right\}.
\end{eqnarray*}
Under ${\cal H}_{0}$, $LR$ converges  to a chi-square distribution with $q$ degrees of freedom,
where $q$ in the number of restrictions imposed by ${\cal H}_{0}$. This approximation can be improved 
if one applies a suitable adjustment to the test statistic. Skovgaard~(2001) proposed two adjusted likelihood ratio
statistics that are asymptotically equivalent for testing ${\cal H}_{0}$. We shall denote them by
$LR^*$ and $LR^{**}$. The adjustment terms depend on a suitable ancillary statistic and involves derivatives with respect to the sample space. A statistic ${\bfsym a}$ is said to be an ancillary statistic if it is distribution constant and, when coupled with the maximum likelihood estimator $\bfsym{\widehat{\theta}}$, is a minimal sufficient statistic for
the model (Barndorff--Nielsen~(1986)). If $(\bfsym{\widehat{\theta}}, {\bfsym a})$ is sufficient, 
but not minimal sufficient, Skovgaard's results still hold; see Severini (2000, Sect. 6.5). In this case,
the log-likelihood function depends on the data only through
$(\bfsym{\widehat{\theta}},{\bfsym a})$ and we write $\ell(\bfsym
{\theta};\bfsym{\widehat \theta},{\bfsym a})$. The sample space derivatives involved are 
$\ell' = \partial \ell(\bfsym{\theta};\bfsym{\widehat \theta},a)/\partial \bfsym{\widehat \theta}$ and
$U'= \partial^2 \ell(\bfsym{\theta};\bfsym{\widehat \theta},
{\bfsym a})/\partial \bfsym{\widehat \theta} \partial{\bfsym{\theta}}^\top$.
The adjusted statistics are given by 
\begin{eqnarray}\label{E.4.15}
LR^* = LR \left(1 - \frac{1}{LR} \log\rho \right)^2
\end{eqnarray}
and
\begin{eqnarray}\label{E.4.16}
LR^{**} = LR - 2\log\rho,
\end{eqnarray}
with 
\begin{eqnarray}\label{E.4.17}
\rho = |\widehat{J}\:|^{1/2} |{\widetilde U}'|^{-1} |{\widetilde J}_{\bfsym{\omega\omega}}|^{1/2} 
|\:{\J}_{\bfsym{\omega\omega}}|^{-1/2} |\:{\J}\:|^{1/2} \frac{\{{\widetilde U}^{\top} {\J}^{\: -1} 
{\widetilde U}\}^{p/2}}{LR^{q/2 - 1} ({\widehat \ell}'- {\widetilde \ell}')^{\top} 
({\widetilde U}')^{-1} {\widetilde U}}.
\end{eqnarray}
Here ${\:\:{\J}}$ equals $\partial^2 \ell(\bfsym{\theta};\bfsym{\widehat \theta},{\bfsym a})/\partial \bfsym{\widehat \theta} \partial{\bfsym{\theta}^\top}$ evaluated at $\bfsym{\widehat \theta} = \bfsym{\widetilde \theta}$ and $\bfsym{\theta} = \bfsym{\widetilde \theta}$. Also, ${\:\:{\J}}_{\bfsym{\omega\omega}}$ is the lower right submatrix of ${\:\:{\J}}$ that corresponds to the nuisance parameter $\bfsym{\omega}$. Both statistics have an approximate 
${\cal X}^2_q$ distribution with high degree of accuracy under the null hypothesis
(Skovgaard, 2001, p. 7). 

Let 
${\bfsym a} = {\bfsym a}({\bfsym z}) = \Big({\bfsym a}_{(1)}^\top({\bfsym z}), \ldots,$ 
${\bfsym a}_{(p)}^\top({\bfsym z})\Big)^\top$, where ${\bfsym a}_{(k)}^\top 
({\bfsym z}) = \left({\bfsym a}_{1k}^\top ({\bfsym z}), \ldots, {\bfsym a}_{{n_k}k}^\top 
({\bfsym z})\right)^\top$, with
$
{\bfsym a}_{jk} ({\bfsym z}) = {\widehat P}_k^{-1}({\bfsym z})\left({\bfsym z}_{jk} - 
{\bfsym{\widehat \mu}}_k({\bfsym z})\right), \ \ j = 1, 2, \ldots, n_k, \ \ k = 1, 2, \ldots, p,
$
where $P_k$ is a lower triangular matrix such that $P_k P_k^\top = \Sigma_k$ is the 
Cholesky decomposition. Following Melo \& Ferrari~(2010) it can be shown that 
${\bfsym a}$ is an ancillary statistic.
With this ancillary statistic we can obtain the sample space derivatives which are
required for the computation of the adjustment term $\rho$. 

In the following we present  some  matrices and vectors that form the score $U$ (Appendix A), the observed information matrix 
$J$ (Appendix A) and the sample space derivatives $\ell'$, $U'$ and ${\:\:{\J}}$ (Appendix B). 
In matrix notation, we have $U = (U_{(1)}^\top, \ldots, U_{(p)}^\top)^\top$, $J = \diag\{J_{(1)}, 
\ldots, J_{(p)}\}$, ${\ell'} = ({\ell}_{(1)}'^\top, \ldots, {\ell}_{(p)}'^\top)^\top$, 
$U' = \diag\{{U}_{(1)}', \ldots, {U}_{(p)}'\}$ and ${\:\:{\J}} = \diag\{{\:\:{\J}}_{(1)}, 
\ldots, {\:\:{\J}}_{(p)}\}$, with $U_{(k)} = - \frac{n_k}{2}\: {\bfsym n}_{(k)}^* + R_{(k)}^\top 
{\bfsym h}_{(k)}, \:\:\: J_{(k)} = \frac{n_k}{2}\: T_{(k)} - R_{(k)}^\top M_{(k)} - V_{(k)}^\top 
Q_{(k)}$, ${\ell'}_{(k)} = 2 R_{(k)}^\top {\bfsym w}_{(k)}, \:\:\: {U'}_{(k)} = 
2(R_{(k)}^\top B_{(k)} + V_{(k)}^\top C_{(k)})$ and ${\:\:{\J}}_{(k)} = 2({\widehat R}_{(k)}^\top 
F_{(k)} + {\widehat V}_{(k)}^\top G_{(k)})$. The $i$-th element of ${\bfsym n}_{(k)}^*$ 
is $\tr({\Sigma}_k^{-1} {\Sigma}_{(k)i})$, for $i = 1, 2, \ldots, s$ and $k = 1, 2, \ldots, p$. 
Here, $s$ is the total number of parameters in ${\bfsym \theta}_{(k)}$. When the ratio
$\lambda_{x_k}$ or the intercept $\alpha_k$ is known, we have $s = 2l+3$, and when the ratio
${\lambda}_{e_k}$ is known, $s = l+4$. The $(i,i')$-th element of $T_{(k)}$ is 
\begin{equation*}\label{E.5.17}
t_{(k)ii'} = \tr({\Sigma}^{(k)i} {\Sigma}_{(k)i'}) + \tr({\Sigma_k}^{-1} \Sigma_{(k)ii'}),
\end{equation*}  
where ${\Sigma}_{(k)i} = {\partial\Sigma_k}/{\partial \theta_{(k)i}}$, $\Sigma_{(k)ii'} = 
{\partial{\Sigma}_{(k)i}}/{\partial \theta_{(k)i'}}$ and ${\Sigma}^{(k)i} = {\partial\Sigma_k^{-1}}/
{\partial \theta_{(k)i}} = - \Sigma_k^{-1} {\Sigma}_{(k)i} \Sigma_k^{-1}$, for $i, i' = 1, 2, 
\ldots, s$. Here, $\theta_{(k)i}$ is the $i$-th element of ${\bfsym \theta}_k$; see 
Appendix C. Also, $R_{(k)}$ and $V_{(k)}$ are block-diagonal matrices given by 
$R_{(k)} = {\rm diag}({\bfsym r}_{(k)}, {\bfsym r}_{(k)}, \ldots, {\bfsym r}_{(k)})$ and
$V_{(k)} = {\rm diag}({\bfsym v}_{(k)}, {\bfsym v}_{(k)}, \ldots, {\bfsym v}_{(k)})$, 
with dimension $sn_k \times s$, and $j$-th element of the vectors ${\bfsym r}_{(k)}$ and
${\bfsym v}_{(k)}$ given by ${r}_{(k)j} = W_{p_0}({\bfsym d}_{jk}^\top 
\Sigma_k^{-1} {\bfsym d}_{jk})$ and ${v}_{(k)j} = W_{p_0}'({\bfsym d}_{jk}^\top \Sigma_k^{-1} 
{\bfsym d}_{jk})$, respectively. Additionally, we define the column vectors ${\bfsym h}_{(k)} = 
\left({\bfsym h}_{(k)}^{(1)^\top}, \ldots, {\bfsym h}_{(k)}^{(s)^\top} \right)^\top$ and 
${\bfsym w}_{(k)} = \left({\bfsym w}_{(k)}^{(1)^\top}, \ldots, {\bfsym w}_{(k)}^{(s)^\top}
\right)^\top$, with dimension $sn_k$, and $j$-th element of the vectors
${\bfsym h}_{(k)}^{(i)}$ and ${\bfsym w}_{(k)}^{(i)}$, for $i = 1, 2, \ldots, s$, given,
respectively, by 
\begin{equation*}\label{E.5.18}
{h}_{(k)j}^{(i)} = {\bfsym d}_{jk}^\top {\Sigma}^{(k)i} {\bfsym d}_{jk} - 
2{\bfsym \mu}_{(k)i}^\top  {\Sigma}_k^{-1} {\bfsym d}_{jk}
\end{equation*}
and
\begin{equation*}\label{E.5.19}
{w}_{(k)j}^{(i)} = \left({\widehat P}_{(k)i} {\bfsym a}_{jk} + 
\bfsym{\widehat \mu}_{(k)i}\right)^\top \Sigma_k^{-1} \left({\widehat P}_k {\bfsym a}_{jk} + 
\bfsym{\widehat \mu}_k - {\bfsym \mu}_k\right),
\end{equation*} 
where ${\widehat P}_{(k)i} = {\partial {\widehat P}_k}/{\partial {\widehat \theta}_{(k)i}}$ and
$\bfsym{\widehat \mu}_{(k)i} = {\partial\bfsym{\widehat \mu}_k}/{\partial {\widehat \theta}_{(k)i}}$. 
The derivative ${\widehat P}_{(k)i}$ is obtained through the algorithm proposed by Smith~(1995) 
and the derivative $\bfsym{\widehat \mu}_{(k)i}$ is presented in Appendix C. The block matrices
$B_{(k)}$, $C_{(k)}$, $F_{(k)}$, $G_{(k)}$, $M_{(k)}$ and $Q_{(k)}$, with dimension $sn_k \times s$, 
have the $(i,i')$-th block given, respectively, by the vectors ${\bfsym b}_{(k)}^{ii'}$, 
${\bfsym c}_{(k)}^{ii'}$, ${\bfsym f}_{(k)}^{ii'}$, ${\bfsym g}_{(k)}^{ii'}$, 
${\bfsym m}_{(k)}^{ii'}$ and ${\bfsym q}_{(k)}^{ii'}$. The $j$-th elements of these vectors 
are, respectively,
\begin{equation*}\label{E.5.20}
\begin{split}
{b}_{(k)j}^{ii'} &= ({\widehat P}_{(k)i'} \ {\bfsym a}_{jk} + \bfsym{\widehat \mu}_{(k)i'})^\top 
{\Sigma}^{(k)i} ({\widehat P}_k {\bfsym a}_{jk} + \bfsym{\widehat \mu}_k - \bfsym{\mu}_k) - 
 \bfsym{\mu}_{(k)i}^\top \Sigma_k^{-1} ({\widehat P}_{(k)i'} \ {\bfsym a}_{jk} + 
 \bfsym{\widehat \mu}_{(k)i'}),
 \\
{c}_{(k)j}^{ii'} &= ({\widehat P}_{(k)i'} \ {\bfsym a}_{jk} + \bfsym{\widehat \mu}_{(k)i'})^\top 
\Sigma_k^{-1} ({\widehat P}_k {\bfsym a}_{jk} + \bfsym{\widehat \mu}_k - \bfsym{\mu}_k) 
\Big[({\widehat P}_k {\bfsym a}_{jk} + \bfsym{\widehat \mu}_k - \bfsym{\mu}_k)^\top {\Sigma}^{(k)i} 
({\widehat P}_k {\bfsym a}_{jk} + \bfsym{\widehat \mu}_k - \bfsym{\mu}_k) 
\\
& \ \ \ \ \ - 2 \bfsym{\mu}_{(k)i}^\top {\Sigma_k}^{-1}({\widehat P}_k {\bfsym a}_{jk} + 
\bfsym{\widehat \mu}_k - \bfsym{\mu}_k) \Big],
\\ \\
{f}_{(k)j}^{ii'} &= ({\widetilde P}_{(k)i'} \ {\bfsym a}_{jk} + \bfsym{\widetilde \mu}_{(k)i'})^\top 
{\widetilde\Sigma}^{(k)i} {\widetilde P}_k {\bfsym a}_{jk} -  \bfsym{\widetilde\mu}_{(k)i}^\top 
{\widetilde\Sigma}_k^{-1} ({\widetilde P}_{(k)i'} \ {\bfsym a}_{jk} + \bfsym{\widetilde \mu}_{(k)i'}),
\end{split}
\end{equation*}

\begin{equation*}\label{E.5.20.1}
\begin{split}
\\
{g}_{(k)j}^{ii'} &= ({\widetilde P}_{(k)i'} \ {\bfsym a}_{jk} + \bfsym{\widetilde \mu}_{(k)i'})^\top 
{\widetilde\Sigma}_k^{-1} {\widetilde P}_k {\bfsym a}_{jk} \left({\bfsym a}^\top_{jk}{\widetilde P}^\top_k 
{\widetilde\Sigma}^{(k)i} {\widetilde P}_k {\bfsym a}_{jk} 
- 2 \bfsym{\widetilde\mu}_{(k)i}^\top {\widetilde\Sigma_k}^{-1} {\widetilde P}_k {\bfsym a}_{jk}\right),
\\ \\
{m}_{(k)j}^{ii'} &= {\bfsym d}_{jk}^\top \Sigma^{(k)ii'} {\bfsym d}_{jk} - 
2\bfsym{\mu}_{(k)i}^\top {\Sigma}^{(k)i'} {\bfsym d}_{jk} - 2\bfsym{\mu}_{(k)i'}^\top {\Sigma}^{(k)i} 
{\bfsym d}_{jk} - 2 \bfsym{\mu}_{(k)ii'}^\top  {\Sigma_k}^{-1} {\bfsym d}_{jk} + 
2 \bfsym{\mu}_{(k)i}^\top {\Sigma_k}^{-1} \bfsym{\mu}_{(k)i'},
\\ \\
{q}_i^{(jk)} &= \left({\bfsym d}_{jk}^\top {\Sigma}^{(k)i} {\bfsym d}_{jk} - 2 \bfsym{\mu}_{(k)i}^\top 
{\Sigma_k}^{-1} {\bfsym d}_{jk} \right) \left({\bfsym d}_{jk}^\top {\Sigma}^{(k)i'} {\bfsym d}_{jk} - 
2\bfsym{\mu}_{(k)i'}^\top {\Sigma_k}^{-1} {\bfsym d}_{jk} \right),
\end{split}
\end{equation*}
where $\bfsym{\mu}_{(k)ii'} = {\partial\bfsym{\mu}_{(k)i}}/{\partial \theta_{(k)i'}}$ and $\Sigma^{(k)ii'} = 
{\partial{\Sigma}^{(k)i}}/{\partial \theta_{(k)i'}} = - 2{\Sigma}^{(k)i} {\Sigma}_{(k)i'} \Sigma_k^{-1} 
- \Sigma_k^{-1} \Sigma_{(k)ii'}\Sigma_k^{-1}$; see Appendix C. 

By replacing $\widehat J$, ${\widetilde J}_{\bfsym{\omega\omega}}$, ${\:{\J}}_{\bfsym{\omega\omega}}$, 
${\:{\J}}$, $\widetilde U$, ${\widetilde U}'$, ${\widehat \ell}'- {\widetilde \ell}'$ and the likelihood ratio
statistic $LR$ in (\ref{E.4.17}) we obtain $\rho$, the quantity that is required for computing the adjusted statistic
$LR^{*}$ in (\ref{E.4.15}) and its equivalent version $LR^{**}$ given in (\ref{E.4.16}).
Note that $\rho$ depends on ${\bfsym Z}_{jk}$, $\bfsym{\mu}_k$, $\Sigma_k$, $\Sigma_k^{-1}$, $P_k$ 
and their first and second derivatives with respect to the parameters. It is worth mentioning that  
the distribution of ${\bfsym Z}_{jk}$ is only required for obtaining the matrices $R_{(k)}$ 
and $V_{(k)}$.

As a final remark, we mention the connection between our results and those obtained by 
Melo \& Ferrari~(2010). In their paper, the model under study is the special case of model 
(\ref{E.5.1}) when $l = p = 1$. The authors obtained the Barndorff-Nielsen~(1986) adjustment 
to the signed likelihood ratio statistic for testing hypotheses on a scalar parameter. 
The adjustment term, given in eq. (6) of their paper, can be calculated using the quantities obtained 
in the present paper for the case in which $\psi$ is scalar ($q=1$). 
Therefore, our results enables us to calculate Barndorff-Nielsen's~(1986) adjusted signed likelihood
ratio statistic in model (\ref{E.5.1}). Hence, our results generalize those in 
Melo \& Ferrari~(2010).

\section{Simulation study}
\label{sec4}

In this section we present a Monte Carlo simulation study to evaluate the efficacy of the adjustments
derived in the previous section. The performances of the tests that use the likelihood ratio 
statistic $(LR)$, and the adjusted statistics  ($LR^{*}$ and $LR^{**}$) 
will be compared with respect to the type I error probability. 

The simulations use model $(\ref{E.5.3})$ with $l = 1$ and $p = 5$. Two different distributions for ${\bfsym Z}_{jk}$ are considered, namely a bivariate normal distribution and a bivariate Student-$t$  with 3 degrees of freedom ($\nu = 3$). 
The number of Monte Carlo replications was 10,000, the nominal levels
of the tests are $\gamma = 1\%, 5\%$ and $10\%$ and the sample sizes are $n_1=\ldots=n_p = 10, 20, 30$ and $40$. 
All simulations were performed using the {\tt Ox} matrix programming language; Doornik~(2006). 

We consider tests of  
$
{\mathcal H}_{0}: {\bfsym \psi} = {\bfsym \psi}^{(0)} \:\:\: {\rm versus} \:\:\: {\mathcal H}_{1}: 
{\bfsym \psi} \neq {\bfsym \psi}^{(0)},
$
where ${\bfsym \psi} = (\beta_1, \beta_2, \ldots, \beta_q)^\top$, for $q = 2, 3, 4, 5$. Also, we 
consider ${\bfsym \psi}^{(0)} = {\bfsym 0}$ when $\lambda_{x_k}$ or $\lambda_{e_k}$ is known, and 
${\bfsym \psi}^{(0)} = {\bfsym 1}$ when the intercept is null. The true parameter values are 
$\alpha_1 = \cdots = \alpha_5 = 0.5$, $\sigma^2_{x_1} = \cdots = \sigma^2_{x_5} = 1.5$, $\sigma^2_{u_1} = 
\cdots = \sigma^2_{u_5} = 0.5$, $\sigma^2_{e_1} = \cdots = \sigma^2_{e_5} = 2.0$. For $\lambda_{x_k}$ or
$\lambda_{e_k}$ known, we set $\mu_{x_1} = \cdots = \mu_{x_5} = 0.5$, and when the intercept is null
we set $\mu_{x_1} = \cdots = \mu_{x_5} = 5.0$.

Tables 1 and 2 present rejection rates (in percentage) of the three tests for all the scenarios described above. 
We notice that the likelihood ratio test ($LR$) is  liberal when the sample size is small in all the cases considered
here. For instance, when ${\bfsym Z}_{jk}$ is normally distributed, $q = 3$, $\lambda_{e_k}$ 
is known and $n_k=10$, the rejection rates of the test that uses $LR$ are $11.4\%$ $(\gamma=5\%)$ and $19.4\%$
$(\gamma=10\%)$; see Table 1. Under the same scenario, except that ${\bfsym Z}_{jk}$ now follows a Student-$t$ 
distribution, the rejection rates are $10.9\%$ $(\gamma=5\%)$ and $18.6\%$ $(\gamma=10\%)$. 
The adjusted tests ($LR^*$ and $LR^{**}$), on the other hand,  display much
better behavior in all cases: they are much less size distorted than the likelihood ratio test.
For example, in the normal case with $\lambda_{x_k}$ known,  $n_k = 10$ and $\gamma=10\%$, 
the rejection rates are $16.9\%$ $(LR)$, $10.2\%$ $(LR^*)$ and $9.8\%$ $(LR^{**})$ for $q = 2$, 
and $18.6\%$ $(LR)$, $10.2\%$ $(LR^*)$ and $9.5\%$ $(LR^{**})$ for $q = 3$. 
As a second example, we mention the case in which the underlying distribution is normal,
the intercept is null, $q = 5$, $n_k = 10$ and $\gamma=5\%$. The rejection rates are  $9.3\%$ 
$(LR)$, $5.2\%$ $(LR^*)$ and $5.0\%$ $(LR^{**})$. Also, for the normal case with 
$\lambda_{e_k}$ known, $n_k = 20$ and $\gamma=1\%$ the rejection rates are $1.9\%$ $(LR)$, 
$1.1\%$ $(LR^*)$ and $1.0\%$ $(LR^{**})$. It can be noticed that, as the number of parameters under test
($q$) grows, the likelihood ratio test deteriorates
while the behavior of the adjusted tests remains unaltered. See, for example, the figures in Table $1$
relative to the Student-$t$ case with $\lambda_{e_k}$ known, $n_k = 10$ and $\gamma = 10\%$; 
the rejection rates are $18.4\%$ $(q=2)$, $18.6\%$ $(q=3)$, $19.9\%$ $(q=4)$ and $21.3\%$ $(q=5)$ for $LR$, $11.4\%$ $(q=2)$, $10.4\%$ $(q=3)$, $10.8\%$ $(q=4)$ and $11.0\%$ $(q=5)$ for $LR^*$, and $10.9\%$ $(q=2)$, $10.1\%$ $(q=3)$, $10.2\%$ $(q=4)$ and
$10.1\%$ $(q=5)$  for $LR^{**}$.

\vspace{0.8cm}
\centerline{\large [Tables 1 and 2 here]}
\vspace{0.8cm}

Our numerical results confirm that the adjusted tests are much better behaved  than the original 
likelihood ratio test in small samples. For almost all the cases, the test that uses the
$LR^{**}$ displays slightly better performance than its asymptotically equivalent version, $LR^*$.

\section{Concluding remarks}
\label{sec6}

In this paper we dealt with the issue of performing hypothesis testing
in an elliptical multivariate errors-in-variables model when the sample size is small. 
We derived modified likelihood ratio statistics that follow very closely a chi-squared
distribution under the null hypothesis. Our approach is based on Skovgaard's~(2001) proposal,
which requires the identification of a suitable ancillary statistic.  
We obtained the required ancillary and all the needed quantities to explicitly write the correction term.
Our simulation results clearly suggested that the adjustment we derived is able to correct the
liberal behavior of the likelihood ratio test in small samples.

\section*{Acknowledgements}

We gratefully acknowledge financial support from FAPESP and CNPq. 


{\footnotesize
\begin{appendix}

\section*{Appendix A. The observed information matrix}
\label{secA}
{
The first derivative of the log-likelihood function for the $k$-th group, $k = 1, 2, 
\ldots, p$, with respect to the parameters is
\begin{eqnarray*}\label{D.1.2}
\frac{\partial \ell_k(\bfsym{\theta})}{\partial \theta_{(k)i}} = -\frac{n_k}{2} \tr
\left({\Sigma_k}^{-1} \Sigma_{(k)i}\right) + \sum_{i=1}^{n_k} W_{p_0}({\bfsym d}_{jk}^\top 
\Sigma_k^{-1} {\bfsym d}_{jk}) \left({\bfsym d}_{jk}^\top \Sigma^{(k)i} {\bfsym d}_{jk} - 2 
\bfsym{\mu}_{(k)i} {\Sigma_k}^{-1} {\bfsym d}_{jk}\right),
\end{eqnarray*}
for $i = 1, 2, \ldots, s$, where $W_{p_0}(u) = {\partial\log p_0(u)}/{\partial u}$, 
$\bfsym{\mu}_{(k)i} = \partial {\bfsym{\mu}_k}/\partial\theta_{(k)i}$, $\Sigma_{(k)i} = 
\partial\Sigma_k/\partial\theta_{(k)i}$ and $\Sigma^{(k)i} = \partial{\Sigma_k}^{-1}/
\partial\theta_{(k)i} = - {\Sigma_k}^{-1} \Sigma_{(k)i} {\Sigma_k}^{-1}$. 
The $(i,i')$-th element of the observed information matrix for the $k$-th group, $J_{(k)}$, 
is given by $J_{(k)ii'} = -{\partial^2 \ell_k(\bfsym{\theta})}/{\partial \theta_{(k)i} \partial 
\theta_{(k)i'}}$, i.e.
\begin{equation*}\label{D.1.4}
\begin{split}
J_{(k)ii'} &= \frac{n_k}{2} \tr\left({\Sigma}^{{(k)i}} \Sigma_{{(k)i'}}\right) + \frac{n_k}{2} 
\tr\left({\Sigma_k}^{-1} \Sigma_{{(k)ii'}}\right) - \sum_{i=1}^{n_k}\Bigg\{W_{p_0}'
\left({\bfsym d}_{jk}^\top {\Sigma_k}^{-1} {\bfsym d}_{jk}\right)\Big({\bfsym d}_{jk}^\top 
\Sigma^{{(k)i}} {\bfsym d}_{jk} 
\\
&- 2 \bfsym{\mu}_{{(k)i}}^\top  {\Sigma_k}^{-1} {\bfsym d}_{jk} \Big) \Big({\bfsym d}_{jk}^\top 
\Sigma^{{(k)i'}} {\bfsym d}_{jk} - 2\bfsym{\mu}_{{(k)i'}}^\top  {\Sigma_k}^{-1} {\bfsym d}_{jk} 
\Big) + W_{p_0}\left({\bfsym d}_{jk}^\top {\Sigma_k}^{-1} {\bfsym d}_{jk}\right) 
\Big({\bfsym d}_{jk}^\top \Sigma^{{(k)ii'}} {\bfsym d}_{jk}  
\\
& - 2\bfsym{\mu}_{{(k)i}}^\top {\Sigma}^{{(k)i'}} {\bfsym d}_{jk} - 2\bfsym{\mu}_{{(k)i'}}^\top 
{\Sigma}^{{(k)i}} {\bfsym d}_{jk} - 2\bfsym{\mu}_{{(k)ii'}}^\top  {\Sigma_k}^{-1} {\bfsym d}_{jk} 
+ 2 \bfsym{\mu}_{{(k)i}}^\top {\Sigma_k}^{-1} \bfsym{\mu}_{{(k)i'}}\Big)\Bigg\},
\end{split}
\end{equation*}
for $i, i' = 1, 2, \ldots, s$ and $k = 1, 2, \ldots, p$, where $W_{p_0}'(u) = 
{\partial W_{p_0}(u)}/{\partial u}$, $\bfsym{\mu}_{(k)ii'} = {\partial\bfsym{\mu}_{(k)i}}/
{\partial \theta_{(k)i'}}$, $\Sigma_{(k)ii'} = {\partial\Sigma_{(k)i}}/{\partial \theta_{(k)i'}}$ 
and $\Sigma^{(k)ii'} = {\partial\Sigma^{(k)i}}/{\partial \theta_{(k)i'}} = - 2{\Sigma}^{(k)i'} 
{\Sigma}_{(k)i} {\Sigma_k}^{-1} - {\Sigma_k}^{-1} \Sigma_{(k)ii'}{\Sigma_k}^{-1}$. 
}

\section*{Appendix B. Sample space derivatives ($\ell'$, $U'$ and ${\:{\J}}$)}
\label{secB}
{
Let ${\bfsym a}$ be the ancillary statistic defined in Section
\ref{sec3} and let us write ${\bfsym z}_{jk} = {\widehat P}_k {\bfsym a}_{jk} 
+ \bfsym{\widehat \mu}_k$. Inserting ${\bfsym z}_{jk}$ in the log-likelihood function  
we have
\begin{eqnarray*}\label{D.2.1}
\ell_k(\bfsym{\theta};\bfsym{\widehat \theta},{\bfsym a}) = -\frac{n_k}{2} \log |\Sigma_k| + 
\sum_{i=1}^{n_k} \log p_0\left(({\widehat P}_k {\bfsym a}_{jk} + \bfsym{\widehat \mu}_k - 
\bfsym{\mu}_k)^\top {\Sigma_k}^{-1} ({\widehat P}_k {\bfsym a}_{jk} + \bfsym{\widehat \mu}_k 
- \bfsym{\mu}_k)\right).
\end{eqnarray*}
Hence,
$\ell' = {\partial \ell(\bfsym{\theta};\bfsym{\widehat \theta},a)}/{\partial \bfsym{\widehat \theta}} 
= \Big({\ell}_{(1)}'^\top, {\ell}_{(2)}'^\top, \ldots, {\ell}_{(p)}'^\top\Big)^\top$, 
where the $i$-th element of the vector ${\ell}_{(k)}'$ is
\begin{equation*}\label{D.2.5}
\begin{split}
{\ell}_{(k)i}' &= 2\:\sum_{i=1}^{n_k} W_{p_0}\left(({\widehat P}_k {\bfsym a}_{jk} + 
\bfsym{\widehat \mu}_{k} - \bfsym{\mu}_{k})^\top {\Sigma_k}^{-1} ({\widehat P}_k {\bfsym a}_{jk} + 
\bfsym{\widehat \mu}_k - \bfsym{\mu}_k)\right) ({\widehat P}_{(k)i} {\bfsym a}_{jk} + 
\bfsym{\widehat \mu}_{(k)i})^\top {\Sigma_k}^{-1} 
({\widehat P}_k {\bfsym a}_{jk} + \bfsym{\widehat \mu}_{k} - \bfsym{\mu}_{k}),
\end{split}
\end{equation*}
for $i = 1, 2, \ldots, s$ and $k = 1, 2, \ldots, p$. Also, we have that $U' = {\partial^2 
\ell(\bfsym{\theta};\bfsym{\widehat \theta},{\bfsym a})}/{\partial \bfsym{\widehat \theta} 
\partial{\bfsym{\theta}^\top}} = \diag\Big\{{U}_{(1)}', {U}_{(2)}', \break \ldots, {U}_{(p)}'\Big\}$, 
where the $(i,i')$-th element of the matrix ${U'}_{(k)}$ is given by 
\begin{equation*}\label{D.2.6}
\begin{split}
{U'}_{(k)ii'}&= 2\sum_{i=1}^{n_k}\Bigg\{W_{p_0}\left({\widehat P}_k {\bfsym a}_{jk} + 
\bfsym{\widehat \mu}_{k} - \bfsym{\mu}_{k})^\top {\Sigma_k}^{-1}({\widehat P}_k {\bfsym a}_{jk} + 
\bfsym{\widehat \mu}_{k} - \bfsym{\mu}_{k})\right)\Big(({\widehat P}_{{(k)i'}} {\bfsym a}_{jk} + 
\bfsym{\widehat \mu}_{{(k)i'}})^\top \Sigma^{{(k)i}} ({\widehat P}_k {\bfsym a}_{jk} + 
\bfsym{\widehat \mu}_k  
\\
&\:\:\:\:- \bfsym{\mu}_k) - {\bfsym{\mu}}_{{(k)i}}^\top {\Sigma_k}^{-1} ({\widehat P}_{{(k)i'}} 
{\bfsym a}_{jk} + \bfsym{\widehat \mu}_{{(k)i'}})\Big) 
+ W'_{p_0}\left(({\widehat P}_k {\bfsym a}_{jk} + \bfsym{\widehat \mu}_{k} - \bfsym{\mu}_{k})^
\top {\Sigma_k}^{-1} ({\widehat P}_k {\bfsym a}_{jk} + \bfsym{\widehat \mu}_{k} - \bfsym{\mu}_{k})\right) 
\\
&\:\:\:\:({\widehat P}_{{(k)i'}} {\bfsym a}_{jk} + \bfsym{\widehat \mu}_{(k)i'})^\top 
{\Sigma_k}^{-1} ({\widehat P}_k {\bfsym a}_{jk} + \bfsym{\widehat \mu}_{k} - \bfsym{\mu}_{k})
\Big(({\widehat P}_k {\bfsym a}_{jk} + \bfsym{\widehat \mu}_{k} - \bfsym{\mu}_{k})^\top 
\Sigma^{{(k)i}}({\widehat P}_k {\bfsym a}_{jk} + \bfsym{\widehat \mu}_{k} - \bfsym{\mu}_{k}) - 
2 \bfsym{\mu}_{{(k)i}}^\top 
\\
&\:\:\:\:{\Sigma_k}^{-1}({\widehat P}_k {\bfsym a}_{jk} + \bfsym{\widehat \mu}_{k} - 
\bfsym{\mu}_{k})\Big)\Bigg\},
\end{split}
\end{equation*}
where ${\widehat P}_{(k)i} = {\partial {\widehat P}_k}/{{\widehat \theta_{(k)i}}}$, 
for $i = 1, 2, \ldots, s$ and $k = 1, 2, \ldots, p$. We also have that
$${\:\:{\J}}= \frac{\partial^2 \ell(\bfsym{\theta};\bfsym{\widehat \theta},{\bfsym a})}
{\partial \bfsym{\widehat \theta} \partial{\bfsym{\theta}^\top}}\Bigg|_{\bfsym{\widehat 
\theta} = \bfsym{\widetilde \theta},\:\: \bfsym{\theta} = \bfsym{\widetilde \theta}} = 
\diag\Big\{{\:\:{\J}}_{(1)}, {\:\:{\J}}_{(2)}, \ldots, {\:\:{\J}}_{(p)}\Big\},$$
where the $(i,i')$-th element of ${\:\:{\J}}_{(k)}$ is given by
\begin{equation*}\label{D.2.7}
\begin{split}
{\:\:{\J}}_{(k)ii'}&= 2\sum_{i=1}^{n_k}\Bigg\{W_{p_0}\left({\bfsym{\widehat d}}_{jk}^\top 
{\widehat \Sigma}_k^{-1} {\bfsym{\widehat d}}_{jk}\right)\left(({\widetilde P}_{(k)i'} 
{\bfsym a}_{jk} +  \bfsym{\widetilde \mu}_{(k)i'})^\top {\widetilde\Sigma}^{(k)i}
{\widetilde P}_k {\bfsym a}_{jk} - {\bfsym{\widetilde\mu}}_{(k)i}^\top {\widetilde\Sigma}_k^{-1} 
({\widetilde P}_{(k)i'} {\bfsym a}_{jk} + \bfsym{\widetilde \mu}_{(k)i'})\right) 
\\
&+ W'_{p_0}\left({\bfsym{\widehat d}}_{jk}^\top {\widehat \Sigma}_k^{-1} 
{\bfsym{\widehat d}}_{jk}\right) 
({\widetilde P}_{(k)i'} {\bfsym a}_{jk} + \bfsym{\widetilde \mu}_{(k)i'})^\top 
{\widetilde\Sigma}_k^{-1} {\widetilde P}_k {\bfsym a}_{jk}\left({\bfsym a}_{jk}^\top 
{\widetilde P}_k^\top {\widetilde\Sigma}^{(k)i}{\widetilde P}_k {\bfsym a}_{jk} - 2 
\bfsym{\widetilde\mu}_{(k)i}^\top {\widetilde\Sigma}_k^{-1}{\widetilde P}_k 
{\bfsym a}_{jk}\right)\Bigg\}.
\end{split}
\end{equation*}
In matrix notation we have
\begin{equation*}\label{D.2.9}
{\ell'}_{(k)} = 2 R_{(k)}^\top {\bfsym w}_{(k)}, \:\:\: {U'}_{(k)} = 2\left(R_{(k)}^\top 
B_{(k)} + V_{(k)}^\top C_{(k)}\right)\:\:\:{\rm and}\:\:\: {\:\:{\J}}_{(k)} = 
2\left({\widehat R}_{(k)}^\top F_{(k)} + {\widehat V}_{(k)}^\top G_{(k)}\right);
\end{equation*}
the elements of the matrices $B_{(k)}$, $C_{(k)}$, $F_{(k)}$, $G_{(k)}$ and of the vector
${\bfsym w}_{(k)}$ are defined in Section \ref{sec3} for $k = 1, 2, \ldots, p$.
}

\section*{Appendix C. Derivatives of the vector $\mu_k$ and of the matrix $\Sigma_k$ with respect 
to the parameters}
\label{secC}
{
When the ratio $\lambda_{x_k} = \sigma^2_{x_k}/\sigma^2_{u_k}$ is known the first derivative of $\bfsym{\mu}_{k}$ 
has elements 
\begin{eqnarray}\label{D.3.1.3}
\bfsym{\mu}_{(k)i} = 
\left\{
     \begin{array}{ll}
         \big(0, \ldots, 0, \underbrace{\mu_{x_k}}_{\mbox{position}\: i}, 0, \ldots, 0)^\top, & 
         \mbox{if}\:\: i = 1, 2, \ldots, l \\ 
         \big({\bfsym 1}^\top, 0\big)^\top, & \mbox{if}\:\: i = l+1 \\  
         \big({\bfsym \beta}_k^\top, 1\big)^\top, & \mbox{if}\:\: i = l+2 \\  
         \:\:\:\:\:\:{\bfsym 0}, & \mbox{if}\:\: i = l+3, l+4, \ldots, s.  
      \end{array}
\right.
\end{eqnarray}
The first derivative of $\Sigma_k$ with respect to the parameter vector ${\bfsym \theta}_{(k)}$ 
is now given for $i = 1, 2, \ldots, s$:
\begin{itemize}
  \item if $i = l+1$ or $i = l+2$, then $\Sigma_{(k)i}$ is null;
	\item if $i = l+3$, we have
	$
         \Sigma_{(k)i} = \left(\begin{array}{cc} 
         {\bfsym \beta}_k {\bfsym \beta}_k^\top \lambda_{x_k}  & {\bfsym \beta}_k \lambda_{x_k} \\ 
         \lambda_{x_k} {\bfsym \beta}_k^\top & \lambda_{x_k} + 1
         \end{array}\right);
  $
  \item is $i = l+4, l+5, \ldots, s$, the elements of  $\Sigma_{(k)i}$ are null except for the $(i-l-3,i-l-3)$-th element, which is given by
    $1$.
  
	\item if $i = 1, \ldots, l$, we have
	\begin{eqnarray*}\label{D.3.1.4}
	\Sigma_{(k)i} = 
  \left\{
         \begin{array}{ll}
         \left(\begin{array}{ccccc} 
         2 \beta_{1k}\lambda_{x_k} \sigma^2_{u_k} & \beta_{2k}\lambda_{x_k} \sigma^2_{u_k} & \ldots &                           \beta_{lk}\lambda_{x_k} \sigma^2_{u_k} & \lambda_{x_k} \sigma^2_{u_k}\\
         \beta_{2k}\lambda_{x_k} \sigma^2_{u_k}   & 0 & \cdots & 0 & 0\\
         \vdots & \vdots & \vdots & \vdots& \vdots\\ 
         \beta_{lk}\lambda_{x_k} \sigma^2_{u_k}   & 0 & \cdots & 0 & 0 \\
         \lambda_{x_k} \sigma^2_{u_k}   & 0 & \cdots & 0 & 0 \\
         \end{array}\right), & 
         \mbox{if}\:\: i = 1 \\ 
%
%
%
         \:\:\:\:\:\:\:\:\:\:\:\:\:\:\:\:\:\:\:\:\:\:\:\:\:\:\:\:\:\:\:\:\:\:\:\:\:\:\:\:\:\:\:\:\:\:\:\:\:\:\:\:\:\:\:\:\:\:\vdots & \:\:\:\:\:\:\:\:\:\:\vdots\\ 
         \left(\begin{array}{ccccc} 
          0 & 0 & \ldots & \beta_{1k}\lambda_{x_k} \sigma^2_{u_k} & 0\\
          0 & 0 & \cdots & \beta_{2k}\lambda_{x_k} \sigma^2_{u_k} & 0\\
         \vdots & \vdots & \vdots & \vdots& \vdots\\ 
         \beta_{1k}\lambda_{x_k} \sigma^2_{u_k}& \beta_{2k}\lambda_{x_k} \sigma^2_{u_k} & \cdots &2                          \beta_{lk}\lambda_{x_k} \sigma^2_{u_k}& \lambda_{x_k} \sigma^2_{u_k} \\
         0 & 0 & \cdots & \lambda_{x_k} \sigma^2_{u_k} & 0 \\
         \end{array}\right), & 
         \mbox{if}\:\: i = l. \\   
      \end{array}
  \right.
  \end{eqnarray*}

\end{itemize}

The second order derivative of $\bfsym{\mu}_{k}$ is
\begin{eqnarray}\label{D.3.1.6}
\bfsym{\mu}_{(k)ii'} = 
\left\{
     \begin{array}{ll}
         \big(0, \ldots, 0, \underbrace{1}_{\mbox{position}\: i}, 0, \ldots, 0)^\top, & \mbox{if}\:\: i = 1, 2, \ldots, l \:\:\mbox{and}\:\: i' = l+2 \\ \\
         \big(0, \ldots, 0, \underbrace{1}_{\mbox{position}\: i'}, 0, \ldots, 0)^\top, & \mbox{if}\:\: i = l+2 \:\:\mbox{and}\:\: i' = 1, 2, \ldots, l \\  
         \:\:\:\:\:\:\:\:\:\:\:\:\:\:\:\:\:\:\:\:\:\:\:\:\:\:\:{\bfsym 0}, & \mbox{otherwise}.  
      \end{array}
\right.
\end{eqnarray}
For $i, i' = 1, 2, \ldots, s$, the elements of the matrix $\Sigma_{(k)ii'}$ are null except for the following cases:
\vspace{-.3cm} 
\begin{itemize}
	\item if $i = 1, 2, \ldots, l$ and $i' = 1, 2, \ldots, l$, we have that, $i = i'$, the elements of $\Sigma_{(k)ii}$ are null except the $(i,i)$-th element, which is equal to $2\lambda_{x_k} \sigma^2_{u_k}$. When $i \neq i'$, the elements of
	the matrix $\Sigma_{(k)ii'}$ are null, except those in positions $(i,i')$ and $(i',i)$, which are equal to  $\lambda_{x_k} \sigma^2_{u_k}$;
	\vspace{-.2cm} 
	\item if $i = 1, \ldots, l$ and $i' = l+3$ we have
	\begin{eqnarray*}\label{D.3.1.7}
	\Sigma_{(k)ii'} = \Sigma_{(k)i'i} = 
  \left\{
         \begin{array}{ll}
         \left(\begin{array}{ccccc} 
         2 \beta_{1k}\lambda_{x_k}&\beta_{2k}\lambda_{x_k} & \ldots & \beta_{lk}\lambda_{x_k} & \lambda_{x_k}\\
         \beta_{2k}\lambda_{x_k} & 0 & \cdots & 0 & 0\\
         \vdots & \vdots & \vdots & \vdots& \vdots\\ 
         \beta_{lk}\lambda_{x_k} & 0 & \cdots & 0 & 0 \\
         \lambda_{x_k} & 0 & \cdots & 0 & 0 \\
         \end{array}\right), & 
         \mbox{if}\:\: i = 1 \\ 
%
%
%
         \:\:\:\:\:\:\:\:\:\:\:\:\:\:\:\:\:\:\:\:\:\:\:\:\:\:\:\:\:\:\:\:\:\:\:\:\:\:\:\:\:\:\:\:\:\vdots & \:\:\:\:\:\:\:\:\:\:\vdots\\ 
         \left(\begin{array}{ccccc} 
          0 & 0 & \ldots & \beta_{1k}\lambda_{x_k}& 0\\
          0 & 0 & \cdots & \beta_{2k}\lambda_{x_k}& 0\\
         \vdots & \vdots & \vdots & \vdots& \vdots\\ 
         \beta_{1k}\lambda_{x_k}& \beta_{2k}\lambda_{x_k}& \cdots &2 \beta_{lk}\lambda_{x_k}& \lambda_{x_k}\\
         0 & 0 & \cdots & \lambda_{x_k}& 0 \\
         \end{array}\right), & 
         \mbox{if}\:\: i = l. \\   
      \end{array}
  \right.
  \end{eqnarray*}
\end{itemize}

When the ratio ${\bfsym \lambda}_{e_k}$ is known the first and second order derivatives of ${{\bfsym \mu}_k}$ are given, respectively, in (\ref{D.3.1.3}) and (\ref{D.3.1.6}), with $s = l+4$. The derivative of $\Sigma_k$ with respect to
the parameter vector ${\bfsym \theta}_{(k)}$ is
\begin{itemize}
	\item if $i = 1, \ldots, l$ we have
	\begin{eqnarray*}\label{D.3.2.2}
	\Sigma_{(k)i} = 
  \left\{
         \begin{array}{ll}
         \left(\begin{array}{ccccc} 
         2 \beta_{1k}\sigma^2_{x_k} & \beta_{2k}\sigma^2_{x_k} & \ldots &\beta_{lk}\sigma^2_{x_k} & \sigma^2_{x_k}\\
         \beta_{2k}\sigma^2_{x_k}   & 0 & \cdots & 0 & 0\\
         \vdots & \vdots & \vdots & \vdots& \vdots\\ 
         \beta_{lk}\sigma^2_{x_k}   & 0 & \cdots & 0 & 0 \\
         \sigma^2_{x_k}   & 0 & \cdots & 0 & 0 \\
         \end{array}\right), & 
         \mbox{if}\:\: i = 1 \\ 
%
%
%
         \:\:\:\:\:\:\:\:\:\:\:\:\:\:\:\:\:\:\:\:\:\:\:\:\:\:\:\:\:\:\:\:\:\:\:\:\:\:\:\:\:\:\:\:\:\:\:\:\:\:\:\:\:\:\:\:\:\:\vdots & \:\:\:\:\:\:\:\:\:\:\vdots\\ 
         \left(\begin{array}{ccccc} 
          0 & 0 & \ldots & \beta_{1k}\sigma^2_{x_k} & 0\\
          0 & 0 & \cdots & \beta_{2k}\sigma^2_{x_k} & 0\\
         \vdots & \vdots & \vdots & \vdots& \vdots\\ 
         \beta_{1k}\sigma^2_{x_k}& \beta_{2k}\sigma^2_{x_k} & \cdots &2\beta_{lk}\sigma^2_{x_k}& \sigma^2_{x_k} \\
         0 & 0 & \cdots & \sigma^2_{x_k} & 0 \\
         \end{array}\right), & 
         \mbox{if}\:\: i = l; \\   
      \end{array}
  \right.
  \end{eqnarray*}
	
	\item if $i = l+1$ or $i = l+2$, the matrix $\Sigma_{(k)i}$ is null;
	\item if $i = l+3$, 
	$
         \Sigma_{(k)i} = \left(\begin{array}{cc} 
         {\bfsym \beta}_k {\bfsym \beta}_k^\top & {\bfsym \beta}_k  \\ 
         {\bfsym \beta}_k^\top & 1
         \end{array}\right);
  $
  \item if $i = l+4 = s$, 
  $
         \Sigma_{(k)i} = \left(\begin{array}{cc} 
         {\bfsym \lambda}_{e_k} & {\bfsym 0} \\ 
         {\bfsym 0} & 1
         \end{array}\right).
  $
\end{itemize}
For $i, i' = 1, 2, \ldots, s$, we have that the elements of $\Sigma_{(k)ii'}$ are null except for the cases:
\vspace{-.3cm} 
\begin{itemize}
	\item if $i = 1, \ldots, l$ and $i' = l+3$, we have
	\begin{eqnarray*}\label{D.3.2.5}
	\Sigma_{(k)ii'} = \Sigma_{(k)i'i} =
  \left\{
         \begin{array}{ll}
         \left(\begin{array}{ccccc} 
         2 \beta_{1k}&\beta_{2k} & \ldots & \beta_{lk} & 1\\
         \beta_{2k} & 0 & \cdots & 0 & 0\\
         \vdots & \vdots & \vdots & \vdots& \vdots\\ 
         \beta_{lk} & 0 & \cdots & 0 & 0 \\
         1 & 0 & \cdots & 0 & 0 \\
         \end{array}\right), & 
         \mbox{if}\:\: i = 1 \\ 
%
%
%
         \:\:\:\:\:\:\:\:\:\:\:\:\:\:\:\:\:\:\:\:\:\:\:\:\:\:\:\:\:\:\vdots & \:\:\:\:\:\:\:\:\:\:\vdots\\ 
         \left(\begin{array}{ccccc} 
          0 & 0 & \ldots & \beta_{1k}& 0\\
          0 & 0 & \cdots & \beta_{2k}& 0\\
         \vdots & \vdots & \vdots & \vdots& \vdots\\ 
         \beta_{1k}& \beta_{2k}& \cdots &2 \beta_{lk}& 1\\
         0 & 0 & \cdots & 1& 0 \\
         \end{array}\right), & 
         \mbox{if}\:\: i = l; \\   
      \end{array}
  \right.
  \end{eqnarray*}
\item if $i = 1, 2, \ldots, l$ and $i' = 1, 2, \ldots, l$, we have that, for $i = i'$, the elements of the matrix $\Sigma_{(k)ii}$ are null except the $(i,i)$-th elements, which is equal to $2\sigma^2_{x_k}$. When $i \neq i'$, the 
elements of the matrix $\Sigma_{(k)ii'}$ are null except the $(i,i')$-th and the  $(i',i)$-th elements, which are equal to $\sigma^2_{x_k}$.
\end{itemize}

When the intercept $\alpha_k$ in known the first order derivative of  $\bfsym{\mu}_{k}$ is
\begin{eqnarray*}\label{D.3.3.2}
\bfsym{\mu}_{(k)i} = 
\left\{
     \begin{array}{ll}
         \big(0, \ldots, 0, \underbrace{\mu_{x_k}}_{\mbox{position}\: i}, 0, \ldots, 0)^\top, & \mbox{if}\:\: i = 1, 2, \ldots, l \\ 
         \big({\bfsym \beta}_k^\top, 1\big)^\top, & \mbox{if}\:\: i = l+1 \\  \\
         \:\:\:\:\:\:{\bfsym 0}, & \mbox{if}\:\: i = l+2, l+3, \ldots, s.  
      \end{array}
\right.
\end{eqnarray*}
The derivative of $\Sigma_k$ with respect to ${\bfsym \theta}_{(k)}$ is now presented for $i = 1, 2, \ldots, s$. We have
\begin{itemize}
  \item if $i = l+1$, the matrix $\Sigma_{(k)i}$ is null;
	\item if $i = l+2$ we have
	\begin{eqnarray*}\label{D.3.3.4}
         \Sigma_{(k)i} = \left(\begin{array}{cc} 
         {\bfsym \beta}_k {\bfsym \beta}_k^\top& {\bfsym \beta}_k\\ 
         {\bfsym \beta}_k^\top &1
         \end{array}\right);
  \end{eqnarray*}
  \item if $i = l+3$ we have
	$
         \Sigma_{(k)i} = \left(\begin{array}{cc} 
         {\bfsym 0}& {\bfsym 0}\\ 
         {\bfsym 0}&1
         \end{array}\right);
  $
  \item if $i = l+4, l+5, \ldots, s$, the elements of the matrix  $\Sigma_{(k)i}$ is null except for the 
   $(i-l-3,i-l-3)$-th element, which equals $1$;
	\item if $i = 1, \ldots, l$, we have
	\begin{eqnarray*}\label{D.3.3.3}
	\Sigma_{(k)i} = 
  \left\{
         \begin{array}{ll}
         \left(\begin{array}{ccccc} 
         2 \beta_{1k}\sigma^2_{x_k} & \beta_{2k}\sigma^2_{x_k} & \ldots &\beta_{lk}\sigma^2_{x_k} & \sigma^2_{x_k}\\
         \beta_{2k}\sigma^2_{x_k}   & 0 & \cdots & 0 & 0\\
         \vdots & \vdots & \vdots & \vdots& \vdots\\ 
         \beta_{lk}\sigma^2_{x_k}   & 0 & \cdots & 0 & 0 \\
         \sigma^2_{x_k}   & 0 & \cdots & 0 & 0 \\
         \end{array}\right), & 
         \mbox{if}\:\: i = 1 \\ 
%
%
%
         \:\:\:\:\:\:\:\:\:\:\:\:\:\:\:\:\:\:\:\:\:\:\:\:\:\:\:\:\:\:\:\:\:\:\:\:\:\:\:\:\:\:\:\:\:\:\:\:\:\:\:\:\:\:\:\:\:\:\vdots & \:\:\:\:\:\:\:\:\:\:\vdots\\ 
         \left(\begin{array}{ccccc} 
          0 & 0 & \ldots & \beta_{1k}\sigma^2_{x_k} & 0\\
          0 & 0 & \cdots & \beta_{2k}\sigma^2_{x_k} & 0\\
         \vdots & \vdots & \vdots & \vdots& \vdots\\ 
         \beta_{1k}\sigma^2_{x_k}& \beta_{2k}\sigma^2_{x_k} & \cdots &2\beta_{lk}\sigma^2_{x_k}& \sigma^2_{x_k} \\
         0 & 0 & \cdots & \sigma^2_{x_k} & 0 \\
         \end{array}\right), & 
         \mbox{if}\:\: i = l. \\   
      \end{array}
  \right.
  \end{eqnarray*}
\end{itemize}

The second order derivative of $\bfsym{\mu}_{k}$ is
\begin{eqnarray*}\label{D.3.3.5}
\bfsym{\mu}_{(k)ii'} = 
\left\{
     \begin{array}{ll}
         \big(0, \ldots, 0, \underbrace{1}_{\mbox{position}\: i}, 0, \ldots, 0)^\top, & \mbox{if}\:\: i = 1, 2, \ldots, l \:\:\mbox{and}\:\: i' = l+1 \\ \\
         \big(0, \ldots, 0, \underbrace{1}_{\mbox{position}\: i'}, 0, \ldots, 0)^\top, & \mbox{if}\:\: i = l+1 \:\:\mbox{and}\:\: i' = 1, 2, \ldots, l \\  
         \:\:\:\:\:\:\:\:\:\:\:\:\:\:\:\:\:\:\:\:\:\:\:\:\:\:\:{\bfsym 0}, & \mbox{otherwise}.  
      \end{array}
\right.
\end{eqnarray*}
For $i, i' = 1, 2, \ldots, s$, we have that the elements of $\Sigma_{(k)ii'}$ are null except for the cases:
\vspace{-.3cm} 
\begin{itemize}
	\item if $i = 1, 2, \ldots, l$ and $i' = 1, 2, \ldots, l$, we have that, for $i = i'$, the elements of $\Sigma_{(k)ii}$
	are null except for the $(i,i)$-th element, which is equal to $2\:\sigma^2_{x_k}$. When $i \neq i'$, the elements of $\Sigma_{(k)ii'}$ are null except for those in positions $(i,i')$ and $(i',i)$, which are equal to $\sigma^2_{x_k}$;
	\vspace{-.2cm} 
	\item{ if $i = 1, \ldots, l$ and $i' = l+2$, we have
	\begin{eqnarray*}\label{D.3.3.6}
	\Sigma_{(k)ii'} = \Sigma_{(k)i'i} =
  \left\{
         \begin{array}{ll}
         \left(\begin{array}{ccccc} 
         2 \beta_{1k}&\beta_{2k} & \ldots & \beta_{lk} & 1\\
         \beta_{2k} & 0 & \cdots & 0 & 0\\
         \vdots & \vdots & \vdots & \vdots& \vdots\\ 
         \beta_{lk} & 0 & \cdots & 0 & 0 \\
         1 & 0 & \cdots & 0 & 0 \\
         \end{array}\right), & 
         \mbox{if}\:\: i = 1 \\ 
%
%
%
         \:\:\:\:\:\:\:\:\:\:\:\:\:\:\:\:\:\:\:\:\:\:\:\:\:\:\:\:\:\:\vdots & \:\:\:\:\:\:\:\:\:\:\vdots\\ 
         \left(\begin{array}{ccccc} 
          0 & 0 & \ldots & \beta_{1k}& 0\\
          0 & 0 & \cdots & \beta_{2k}& 0\\
         \vdots & \vdots & \vdots & \vdots& \vdots\\ 
         \beta_{1k}& \beta_{2k}& \cdots &2 \beta_{lk}& 1\\
         0 & 0 & \cdots & 1& 0 \\
         \end{array}\right), & 
         \mbox{if}\:\: i = l. \\   
      \end{array}
  \right.
  \end{eqnarray*}
  }
\end{itemize}
}
\end{appendix}

}

\newpage

\begin{table}[h]
\hspace{1.6cm} {\caption{\normalsize Null rejection rates; $n_k = 10$.}} 
\vspace{0.5cm} \centering
\begin{tabular}{cccccccccccccccc} 
\hline\hline 
\multicolumn{16}{c}{$\lambda_{x_k}$ known}\\\hline 
\multicolumn{8}{c}{Normal distribution}&\multicolumn{8}{c}{Student-$t$ distribution $(\nu = 3)$}\\\hline 
&\multicolumn{3}{c}{$\gamma = 5\%$}&&\multicolumn{3}{c}{$\gamma = 10\%$}&&\multicolumn{3}{c}{$\gamma = 
5\%$}&&\multicolumn{3}{c}{$\gamma = 10\%$}
\\\cline{2-4}\cline{6-8}\cline{10-12}\cline{14-16}
{$q$}& $LR$  & $LR^*$& $LR^{**}$ && $LR$  & $LR^*$& $LR^{**}$   && $LR$  & $LR^*$& $LR^{**}$ && $LR$  
& $LR^*$& $LR^{**}$\\\hline 
2 &10.1&5.1&4.9&&16.9&10.2& 9.8  && 9.9&5.2&5.0&&16.7&10.6&10.2\\
3 &11.6&5.3&5.0&&18.6&10.2& 9.5  &&10.7&5.2&4.9&&18.4&10.3& 9.8\\
4 &12.5&5.3&4.8&&20.4&10.2& 9.6  &&12.0&5.4&5.1&&19.9&10.6&10.1\\
5 &13.6&5.2&4.8&&21.5&10.3& 9.6  &&12.7&5.5&5.1&&20.8&10.4& 9.8\\
\hline\hline
\multicolumn{16}{c}{$\lambda_{e_k}$ known}\\\hline 
\multicolumn{8}{c}{Normal distribution}&\multicolumn{8}{c}{Student-$t$ distribution $(\nu = 3)$}\\\hline 
&\multicolumn{3}{c}{$\gamma = 5\%$}&&\multicolumn{3}{c}{$\gamma = 10\%$}&&\multicolumn{3}{c}{$\gamma = 
5\%$}&&\multicolumn{3}{c}{$\gamma = 10\%$}
\\\cline{2-4}\cline{6-8}\cline{10-12}\cline{14-16}
{$q$}& $LR$  & $LR^*$& $LR^{**}$ && $LR$  & $LR^*$& $LR^{**}$   && $LR$  & $LR^*$& $LR^{**}$ && $LR$  
& $LR^*$& $LR^{**}$\\\hline 
2 &10.1&5.1&4.9&&17.3&10.1& 9.6  &&10.7&5.8&5.5&&18.4&11.4&10.9\\
3 &11.4&5.3&4.9&&19.4&10.1& 9.6  &&10.9&5.3&5.0&&18.6&10.4&10.1\\
4 &13.0&5.2&4.7&&20.8&10.7& 9.8  &&11.9&5.5&5.2&&19.9&10.8&10.2\\
5 &13.7&5.2&4.7&&22.4&10.3& 9.6  &&13.2&5.5&5.2&&21.3&11.0&10.1\\
\hline\hline
\multicolumn{16}{c}{null intercept}\\\hline 
\multicolumn{8}{c}{Normal distribution}&\multicolumn{8}{c}{Student-$t$ distribution $(\nu = 3)$}\\\hline 
&\multicolumn{3}{c}{$\gamma = 5\%$}&&\multicolumn{3}{c}{$\gamma = 10\%$}&&\multicolumn{3}{c}{$\gamma = 
5\%$}&&\multicolumn{3}{c}{$\gamma = 10\%$}
\\\cline{2-4}\cline{6-8}\cline{10-12}\cline{14-16}
{$q$}& $LR$  & $LR^*$& $LR^{**}$ && $LR$  & $LR^*$& $LR^{**}$   && $LR$  & $LR^*$& $LR^{**}$ && $LR$  
& $LR^*$& $LR^{**}$\\\hline 
2 &7.8&5.1&4.9&&13.9&10.2&10.0  &&7.3&5.5&5.4&&13.3&10.6&10.5\\
3 &8.6&5.0&4.9&&15.2&10.2&10.0  &&7.4&5.2&5.2&&13.6&10.4&10.2\\
4 &8.9&5.1&5.0&&15.7&10.2&10.0  &&7.3&4.8&4.7&&13.9&10.2&10.0\\
5 &9.3&5.2&5.0&&16.4&10.1& 9.8  &&7.9&5.1&5.0&&14.6&10.3&10.2\\
\hline
\end{tabular}
\label{tab.1}
\end{table}

\begin{table}[h]
\hspace{1.6cm} {\caption{\small Null rejection rates; $q = 3$. Normal distribution.}} 
\vspace{0.5cm} \centering
\begin{tabular}{cccccccccccc} 
\hline\hline 
\multicolumn{12}{c}{$\lambda_{x_k}$ known}\\
\hline
&\multicolumn{3}{c}{$\gamma = 1\%$}&&\multicolumn{3}{c}{$\gamma = 5\%$}&&\multicolumn{3}{c}
{$\gamma = 10\%$}
\\\cline{2-4}\cline{6-8}\cline{10-12}
{$n_k$}& $LR$ & $LR^*$& $LR^{**}$ && $LR$  & $LR^*$& $LR^{**}$ && $LR$  & $LR^*$& $LR^{**}$ \\\hline 
  10 &3.6&1.1&1.0&&11.6&5.3&5.0&&18.6&10.2& 9.5 \\
  20 &1.9&0.9&0.9&& 7.8&5.1&5.1&&13.8&10.0& 9.9 \\ 
  30 &1.4&1.0&0.9&& 6.3&4.7&4.7&&12.3& 9.7& 9.6 \\ 
  40 &1.3&0.9&0.9&& 6.0&4.8&4.8&&11.7& 9.9& 9.9 \\   
\hline\hline
\multicolumn{12}{c}{$\lambda_{e_k}$ known}\\
\hline
&\multicolumn{3}{c}{$\gamma = 1\%$}&&\multicolumn{3}{c}{$\gamma = 5\%$}&&\multicolumn{3}{c}{$\gamma = 10\%$}
\\\cline{2-4}\cline{6-8}\cline{10-12}
{$n_k$}& $LR$ & $LR^*$& $LR^{**}$ && $LR$  & $LR^*$& $LR^{**}$ && $LR$  & $LR^*$& $LR^{**}$ \\\hline 
  10 &3.5&1.2&1.0&&11.4&5.3&4.9&&19.4&10.1& 9.6 \\ 
  20 &1.9&1.1&1.0&& 8.0&5.2&5.1&&14.6&10.4&10.3 \\ 
  30 &1.4&0.9&0.9&& 6.9&5.2&5.2&&12.4&10.0&10.0 \\ 
  40 &1.3&0.9&0.9&& 6.0&4.8&4.8&&11.7& 9.9& 9.9 \\   
\hline\hline
\multicolumn{12}{c}{null intercept}\\
\hline
&\multicolumn{3}{c}{$\gamma = 1\%$}&&\multicolumn{3}{c}{$\gamma = 5\%$}&&\multicolumn{3}{c}{$\gamma = 10\%$}
\\\cline{2-4}\cline{6-8}\cline{10-12}
{$n_k$}& $LR$ & $LR^*$& $LR^{**}$ && $LR$  & $LR^*$& $LR^{**}$ && $LR$  & $LR^*$& $LR^{**}$ \\\hline 
  10 &2.1&1.1&1.0&&8.6&5.0&4.9&&15.2&10.2&10.0 \\ 
  20 &1.5&0.9&0.9&&6.0&4.5&4.5&&11.9& 9.8& 9.7 \\ 
  30 &1.3&1.1&1.1&&6.0&5.1&5.1&&11.6&10.0& 9.9 \\ 
  40 &1.4&1.1&1.1&&5.6&4.9&4.9&&11.1&10.1&10.1 \\   
\hline
\end{tabular}
\label{tab.2}
\end{table}

%
%
%
%
%

\end{document}